\documentstyle{mathnach} 
\Year{1995}
\Received{December 7, 1995}           

\newcounter{saveeqn}
\def\thealphequation{\mbox{\arabic{section}.\arabic{saveeqn}\alph{equation}}}

\newcommand{\alpheqn}{\setcounter{saveeqn}{\value{equation}}%
    \setcounter{equation}{0}\stepcounter{saveeqn}%
    \renewcommand{\theequation}{\thealphequation}}
\newcommand{\reseteqn}{\setcounter{equation}{\value{saveeqn}}%
     \renewcommand{\theequation}{\arabic{section}.\arabic{equation}}}

\newcommand{\myref}[1]{{\rm (\ref{#1})}}

\newcommand{\plref}[1]{{\rm \ref{#1}}}

\def\ga{\alpha} 
\def\gb{\beta}  
 
  \def\pl{\partial}

\def\gl{\lambda}

\def\eps{\varepsilon}

\def\cd{{\cal D}}

\def\cf{{\cal F}}

\def\ch{{\cal H}}

\def\arg{{\rm arg\,}}

\def\ker{{\rm ker\,}}

\def\log{{\rm log\,}}

\def\Re{{\rm Re\,}}
\def\Res{{\rm Res}}

\def\spec{{\rm spec\,}}

\def\supp{{\rm supp\,}}

\def\Tr{{\rm Tr}}

\def\max{{\rm max}}
\def\min{{\rm min}}

\newcommand{\cinfz}[1]{C_0^\infty(#1)}
\newcommand{\cinf}[1]{C^\infty(#1)}

\def\moplus{\mathop{\oplus}}

\newcommand{\superwidearray}{\renewcommand{\arraystretch}{2}}
\newcommand{\widearray}{\renewcommand{\arraystretch}{1.5}}
\newcommand{\restarray}{\renewcommand{\arraystretch}{1}}
\newcommand{\casetwo}[4]%
{\widearray
  \textstyle\left\{\begin{array}{ll} {#1,} & {#2}\\{#3,} & {#4}\end{array}\right.
  \restarray}

\newcommand{\smallcasetwo}[4]%
{\textstyle\left\{\begin{array}{ll} {#1,} & {#2} \\{#3,} & {#4}\end{array}\right.}

\def\regint{-\hspace*{-1em}\int}
\def\reginttext{-\hspace*{-0.9em}\int}

\def\DST{\displaystyle}
\def\TST{\textstyle}

\newcommand{\casethree}[6]{\mbox{$\left\{\begin{array}{r@{,\quad}l}
 {#1} & {#2}\\{#3} & {#4}\\{#5} & {#6}\end{array}\right.$}}

\def\LIM{\mathop{{\rm LIM}}}

\newcommand{\detz}{{\rm det}_\zeta}

\newcommand{\mltilde}{\widetilde}
\def\secfive{5}

\begin{document}
\Title{Determinants of regular singular Sturm--Liouville operators}
\Shorttitle{Determinants of Sturm--Liouville operators}
\By{{\sc Matthias Lesch } of Berlin }
\Names{Lesch}
\Dedicatory{ }
\Subjclass{Primary 58G26; Secondary 58G11, 34B24.}
\Keywords{Determinants, Sturm--Liouville theory.}
\Email{lesch@mathematik.hu-berlin.de}
\maketitle
\begin{abstract}
We consider a regular singular Sturm--Liouville operator
$$L:=-\frac{d^2}{dx^2} + \frac{q(x)}{x^2 (1-x)^2}$$
on the line segment $[0,1]$. We impose certain boundary
conditions such that we obtain a semi--bounded self--adjoint operator.
It is known (cf. Theorem 1.1 below)
that the $\zeta$--function of this operator
$$\zeta_L(s)=\sum_{\lambda\in\spec(L)\setminus\{0\}} \lambda^{-s}$$
has a meromorphic continuation to the whole complex plane with
$0$ being a regular point. Then, according to \cite{RS} the $\zeta$--regularized
determinant of $L$ is defined by
$$\detz(L):=\exp(-\zeta_L'(0)).$$
In this paper we are going to express this determinant in terms of
the solutions of the homogeneous differential equation $Ly=0$
generalizing earlier work of
{\sc S. Levit and U. Smilansky} \cite{LS}, {\sc T. Dreyfus and H. Dym}
\cite{DD}, and {\sc D. Burghelea, L. Friedlander and T. Kappeler}
\cite{BFK1,BFK2}. More precisely we prove the formula
$$\detz(L)=\frac{\pi\,W(\psi,\varphi)}%
    {2^{\nu_0+\nu_1} \Gamma(\nu_0+1)\Gamma(\nu_1+1)}.$$
Here $\varphi, \psi$ is a certain 
fundamental system of solutions for the
homogeneous equation $Ly=0$, 
$W(\varphi, \psi)$ denotes their Wronski determinant,
and $\nu_0, \nu_1$ are numbers related to the characteristic roots
of the regular singular points $0, 1$.

\end{abstract}  

\newsection{Introduction and statement of the main result}

We begin with some elementary remarks on $\zeta$--regularized determinants.
Let $L\ge -c+1$ be a semi--bounded self--adjoint operator in the Hilbert space
$\ch$. We assume that
\begin{equation}
  (L+c)^{-1}\in C_1(\ch)
\label{G1.1}
\end{equation}
is trace class. Usually one deals with the more general assumption
that only some power of $(L+c)^{-1}$ is trace class. But since
in this paper we deal exclusively with one--dimensional Sturm--Liouville
operators, we may content ourselves with the more convenient
case \myref{G1.1}. Moreover, we assume that we have an asymptotic expansion
\begin{equation}
    \Tr(L+z^2)^{-1}\sim_{z\to\infty}\sum_{j=0}^\infty\sum_{k=0}^{k_j}
  A_{jk}\, z^{\alpha_j-2}\,\log^k z,
  \label{G1.2}
\end{equation}
as $z\to\infty, z\in Z:=\{z\in \C\,|\, |\arg z| < \delta\},\quad
0<\delta< \pi/2$ fixed. Here $k(j)\in\Z_+$ for all $j$,
$(\alpha_j)_{j\in\Z_+}$ is a sequence of complex numbers with
$\Re(\alpha_j)\to -\infty$ (cf. \cite[Sec. 2]{BL94}). We assume
furthermore, that the terms $z^{-2}\, \log^k z, k\ge 1,$ do not
occur, i.e. $k(j)=0$ if $\alpha_j=0$. These assumptions guarantee, that
the $\zeta$--function of $L$,
\begin{equation}
  \zeta_L(s)=\sum_{\lambda\in \spec L\setminus \{0\}} \lambda^{-s},\quad
   \Re s>>0,
\label{G1.3}
\end{equation}
has a meromorphic continuation to $\C$ with no pole at
$0$ (we put $\gl^{-s}=e^{-i\pi s} |\gl|^{-s}$ if $\gl<0$).

Now we can define the $\zeta$--regularized determinant: put
$\detz L=0$ if $0\in\spec L$, otherwise put
\begin{equation}
  \detz L:= \exp(-\zeta_L'(0)).
  \label{G1.4}
\end{equation}

This notion was introduced by Ray--Singer \cite{RS}.

Since we mostly deal with $\log\detz L$ instead of $\detz L$ we 
abbreviate
\begin{equation}
T(L):= -\zeta_L'(0)
\end{equation}
which is $\log\detz L$ for invertible $L$.

If $L\ge 0$, our assumption \myref{G1.1} implies
\begin{equation}
  \zeta_L(s)=2\frac{\sin \pi s}{\pi} \regint_0^\infty
     z^{1-2s}\;\Tr(L+z^2)^{-1} dz
\label{G1.5}
\end{equation}
from which one easily derives the formula
\begin{equation}
  T(L)=-2\regint_0^\infty z\; \Tr(L+z^2)^{-1} dz.
\label{G1.6}
\end{equation}
The symbol $\reginttext$ indicates that the integral has to be
regularized. For convenience of the reader we briefly recall the
definition of $\reginttext$ as we will make extensive use of this
notion: if $f:\R_+\to\C$ is a function having an asymptotic expansion
\begin{equation}
f(x)\sim_{x\to 0} \sum_{j=1}^N\sum_{k=0}^{k_j} a_{jk} \,x^{\ga_j}\,
\log^k x+ \sum_{k=0}^{k_0} a_{k0}\,\log^k x+ o(1)
\label{G1.30}
\end{equation}
with $\Re \ga_j\not=0$, then we define a "regularized limit" by
\begin{equation}
\LIM_{x\to 0} f(x):= a_{00}.
\end{equation}
If $f(x)$ also has an asymptotic expansion of the form \myref{G1.30}
as $x\to\infty$, then $\displaystyle\LIM_{x\to \infty}\textstyle f(x)$
is defined likewise.

Next let $f:\R_+\to\C$, such that 
\begin{equation}
   \begin{array}{rcl}
   \DST  f(x)&=&\DST \sum_{j=1}^N\sum_{k=0}^{k_j} a_{jk}\,x^{\ga_j}\,\log^k x
                +x^\eps f_1(x),\\
         &=&\DST \sum_{j=1}^M\sum_{l=0}^{l_j} b_{jl}\,x^{\gb_j}\,\log^k x+
              x^{-\eps} f_2(x),
   \end{array}
\end{equation}
with
$f_1\in L^1_{\rm loc}([0,\infty)), f_2\in L^1([1,\infty)), \eps>0$.

Then we put
\begin{equation}
        \regint_0^\infty f(x)dx =\LIM_{b\to \infty}\LIM_{a\to 0}
            \int_a^b f(x) dx .
            \label{G1.31}
\end{equation}
This is the partie--fini definition of Hadamard. $\reginttext_0^\infty f(x) dx$
can also be expressed in terms of the Mellin transform
(cf. \cite{BS85,L93}). Namely,
$$F_1(s)=\int_0^1 x^s f(x) dx,\quad F_2(s)=\int_1^\infty x^s f(x)dx$$
extend to meromorphic functions in a half--plane containing $0$.
Then
\begin{equation}
         \regint_0^\infty f(x) dx=\begin{array}{l}
            \mbox{\rm constant term in the Laurent expansion}\\
            \mbox{\rm of } F_1(s)+F_2(s)\, \mbox{\rm at}\, s=0.
            \end{array}
         \label{G1.32}
\end{equation}

One easily calculates
\alpheqn
\begin{eqnarray}
    \regint_0^1 x^\ga \, \log^k x dx&=&\smallcasetwo{0}{\ga=-1,}{\frac{(-1)^k k!}%
          {(\ga+1)^{k+1}}}{\ga\not=-1,}\label{G4-1.13}       \\
    \regint_1^\infty x^\ga \, \log^k x dx&=&\smallcasetwo{0}{\ga=-1,}{\frac{(-1)^{k+1} k!}%
          {(\ga+1)^{k+1}}}{\ga\not=-1,}
\end{eqnarray}
in particular
\begin{equation}
\regint_0^\infty x^\ga \, \log^k x\, dx=0.
\end{equation}
\reseteqn
For an elaborate discussion of $\reginttext$ see \cite[Sec. 2.1]{L93}.

Next we consider a Sturm--Liouville operator
\begin{equation}
  l=-\frac{d^2}{dx^2}+q(x)
\label{G1.7}
\end{equation}
on the interval $(0,1)$, where $q\in\cinf{0,1}$ is a real function.

Assume for the moment that $q\in\cinf{[0,1]}$ and impose, for
simplicity, Dirichlet boundary conditions. More precisely, we consider
the self--adjoint extension, $L$, of $l$ with domain
$$\cd(L)=\{ f\in H^2[0,1]\,|\, f(0)=f(1)=0\}\subset L^2[0,1].$$
Then $L$ satisfies \myref{G1.1}, \myref{G1.2}, hence its $\zeta$--regularized
determinant is well--defined. According to \cite{BFK2},
$\detz(L)$ can be computed as follows: let $\varphi$ be the unique
function with
$$ l\varphi=0,\quad \varphi(0)=0,\quad \varphi'(0)=1.$$
Then
\begin{equation}
  \detz(L)=2\varphi(1).
\end{equation}

In this paper we want to generalize this result to Sturm--Liouville
operators having regular singularities at $0$ and $1$.

From now on let $l$ be regular singular, i.e.
\begin{equation}
   q(x)=x^{-2} a_0(x^{1/N})=(1-x)^{-2} a_1((1-x)^{1/N})
\label{G1.8}
\end{equation}
with $a_0\in\cinf{[0,1)}, a_1\in\cinf{(0,1])}$ and
\begin{equation}
  a_0(0),\, a_1(1)\ge -1/4
\label{G1.9}
\end{equation}
(cf. \cite[Sec. 4]{BL94} for examples where this kind of potential
occurs naturally). For later purposes we write
\begin{equation}
    a_0(x)=:a_0(0)+ x b_0(x);\quad a_1(x)=:a_1(1) +(1-x) b_1(x)
    \label{G1.15}
\end{equation}
with $b_0\in\cinf{[0,1)}, b_1\in\cinf{(0,1]}$.

The assumptions \myref{G1.9} imply that $l$ is bounded below
on $\cinfz{0,1}\subset L^2[0,1]$.
In the sequel, the domain of an operator is denoted by
$\cd$ and we put 
\begin{equation}\begin{array}{rcl}
   \DST   l_\min&:=&\DST\overline{l}=\,\mbox{\rm closure of $l$ in}\; L^2[0,1],\\
   \DST   l_\max&:=&\DST l^*.\end{array}
\end{equation}

If $a_0(0)
\ge 3/4\; (\mbox{resp.}\, a_1(1)\ge 3/4)$ then $0\,
(\mbox{resp.}\, 1)$ is in the limit point case and no boundary condition
is necessary. Otherwise we have to impose boundary conditions to obtain
a self--adjoint operator. Since we will be dealing exclusively with
separated boundary conditions, it is enough to discuss one boundary
point, e.g. $0$. We distinguish between two cases:

\paragraph{1. $q$ is continuous at $0$:} we impose a boundary condition
at $0$ of the form
\begin{equation}
   R_0f=0\quad\mbox{\rm for}\quad f\in \cd(l_\max),
   \label{G1.10}
\end{equation}
where either
\begin{equation}
    R_0f= f(0)\quad\mbox{\rm (Dirichlet)}
    \label{G1.16}
\end{equation}
or
\begin{equation}
    R_0f= f'(0)+Af(0)\quad\mbox{\rm (generalized Neumann)}
    \label{G1.17}
\end{equation}
with some $A\in\R$.

Furthermore we define the order of the boundary operator $R_0$ by
\begin{equation}
     \sigma(R_0):=\casetwo{0}{\mbox{\rm in case}\, \myref{G1.16},}{%
                  1}{\mbox{\rm in case}\, \myref{G1.17}.}
    \label{G1.18}
\end{equation}

\paragraph{2. $q$ is not continuous at $0$:} in this situation we content
ourselves with the 'Dirichlet' condition at $0$. Since $l$ is bounded below,
we can form its Friedrichs extension, $l^\cf$. Now for $f\in\cd(l_\max)$
we require
\begin{equation}
   \varphi f\in\cd(l^\cf)
   \label{G1.11}
\end{equation}
for any cut--off function $\varphi\in\cinfz{[0,1)}$ with $\varphi=1$
in a neighborhood of $0$. In other words we consider 'the Friedrichs
extension near zero'. It can be checked that for $f\in\cd(l_\max)$
the condition \myref{G1.11} is equivalent to
\begin{equation}
    |f(x)|=O(x^{1/2}),\; x\to 0.
    \label{G1.12}
\end{equation}
However, we would like to express this boundary condition in terms of a boundary
operator. For doing this we put
\begin{equation}
   \nu_0:=\sqrt{a_0(0)+1/4}.
\end{equation}
$1/2\pm \nu_0$ are just the roots of the indicial equation
\begin{equation}
\gl (\gl-1) - a_0(0)=0
\end{equation}
of the regular singular point $0$. It is well--known that
there is a fundamental system $\varphi, \psi$ of solutions of
the homogeneous equation $lf=0$, where
\alpheqn
\begin{eqnarray}
         \varphi(x)&=& x^{\nu_0+1/2}\varphi_1(x^{1/N}),\\
         \psi(x)&=&x^{-\nu_0+1/2}\,\psi_1(x^{1/N})+ k \varphi(x)\,\log(x)
\end{eqnarray}
\reseteqn
with $\varphi_1, \psi_1\in\cinf{[0,1)}, \varphi_1(0)=1,$ and
$$ \psi_1(0)=\smallcasetwo{1/2\nu_0}{\nu_0>0,}{0}{\nu_0=0.}$$
Note that $k=0$ if $2N\nu_0\not\in\Z$. If $\nu_0=0$ then $0$ is a repeated
root of the indicial equation and hence we can choose $k=-1$.
We call such a fundamental system {\it normalized at $0$}.
Obviously, the Wronskian of $\varphi, \psi$, $W(\varphi,\psi)=-1$.
We also introduce
\alpheqn
\begin{eqnarray}
         \varphi_{\nu_0}(x)&=& x^{\nu_0+1/2},\label{G4-1.25}\\
         \psi_{\nu_0}(x)&=&\casetwo{x^{-\nu_0+1/2}/2\nu_0}{\nu_0>0,}{%
           -x^{1/2}\,\log(x)}{\nu_0=0,}
\end{eqnarray}
\reseteqn
which is a fundamental system of solutions of the
differential equation
\begin{equation}
l_{\nu_0}f:=-\frac{d^2}{dx^2}f+x^{-2}(\nu_0^2-1/4)f=0.
\label{G1.14}
\end{equation}
Now we put for $f\in\cd(l_\max)$
\begin{equation}
    R_0f:= \lim_{x\to 0} \psi_{\nu_0}(x)^{-1} f(x).
    \label{G1.19}
\end{equation}
Note that $R_0$ is even well--defined on the larger space
\begin{equation}
   \widetilde \cd(l):=\{ f+ c_1 \varphi + c_2 \psi\,|\, f\in\cd(l_\max),\,
     c_1,c_2\in\C\,\}.
   \label{G1.20}
\end{equation}
Moreover if
$$ f= a\varphi + b \psi$$
then
\begin{equation}
R_0f=\lim_{x\to 0} \psi_{\nu_0}(x)^{-1}f(x)=b.
\end{equation}
Finally we define the 'order' of this boundary operator to be
\begin{equation}
   \sigma(R_0):= 1/2-\nu_0.
   \label{G1.21}
\end{equation}

In order to treat the various boundary conditions in a unified
way, we extend the definition of $\nu_0$ to
continuous $q$ and Neumann boundary conditions. We put
$$\nu_0= \frac 12 - \sigma(R_0).$$
Summing up we have
\begin{equation}\superwidearray\nu_0=\left\{\begin{array}{r@{,\quad}l}
       {-\frac 12}&{\mbox{\rm\parbox{6cm}{ if $q$ is continuous at $0$
        and the boundary condition is of type \myref{G1.17},}}}\\
     {\frac 12}&{\mbox{\rm \parbox{6cm}{if $q$ is continuous at $0$
        and the boundary condition is of type \myref{G1.16},}}}\\
     {\sqrt{a_0(0)+1/4}}&{\mbox{\rm \parbox{6cm}{if $q$ is of type \myref{G1.8} 
      and not continuous at $0$,}}}
     \end{array}\right.
\end{equation}
and
$\sigma(R_0)=\frac 12 - \nu_0$. $\nu_1$ and $R_1$ are defined analogously.
\bigskip

Having chosen boundary operators $R_0, R_1$ of the above types we obtain a
self--adjoint extension, $L:=(l, R_0, R_1)$, of $l$ with domain
\begin{equation}
    \cd(L)=\{ f\in\cd(l_\max)\,|\, R_0 f=0, R_1 f=0\}.
    \label{G1.22}
\end{equation}

Our aim is to compute the $\zeta$--regularized determinant of $L$. The
existence of $\detz L$ is a consequence of the following result:

\pagebreak[3]
\begin{theorem}\label{S1.1} {\rm (\cite{BS85,Gi,BL94})} $L$
is a discrete operator, $(L+z^2)^{-1}$ is trace class and we have an
asymptotic expansion
\begin{eqnarray*}
  \Tr((L+z^2)^{-1})&\sim_{z\to\infty}&
    \sum_{j=0}^\infty a_j\; z^{-1-2j}+\sum_{j=1}^\infty
                  b_j\; z^{-1-2j}\,\log z+\sum_{j=0}^\infty
                      c_j\; z^{-2-j/N}\\
         &\sim_{z\to\infty}& a_0 z^{-1} +c_0 z^{-2} +O(z^{-2-1/N}\,\log z).
\end{eqnarray*}
Moreover
$$a_0=1/2$$
and
$$c_0=(\sigma(R_0)+\sigma(R_1)-1)/2=-\frac 12 (\nu_0+\nu_1).$$
\end{theorem}
For smooth potentials this result is classical (cf. \cite[Sec. 1.7]{Gi}). For regular singular
operators it is due to \cite{BS85} in case of $N=1$. For the
extension to arbitrary $N$ see \cite{BL94}.

\def\trlz{\Tr (L+z^2)^{-1}}
If $L\ge 0$ then Theorem \plref{S1.1} implies in view of \myref{G1.6}
\begin{eqnarray}
     T(L)&=& -2 \regint_0^\infty z \;\big[\Tr(L+z^2)^{-1}-
           a_0 z^{-1}-c_0\, z^{-2}\big] \, dz\nonumber \\
       &=& -2 \regint_0^1 z \;\trlz - a_0 \;dz -2 \int_1^\infty
           z \;\trlz - a_0 -c_0 \,z^{-1}\; dz\nonumber \\
      &=&-2 \int_0^1\big[ z \;\trlz - a_0 -z^{-1}\;\dim\ker L \big] \, dz
          \label{G1.38}\\
      && -2 \int_1^\infty \big[ z\; \trlz - a_0 - c_0\, z^{-1}\big ]\, dz
         \nonumber
\end{eqnarray}

Here, we have used $\reginttext_0^1 z^{-1}\,dz=0, \reginttext_0^\infty
z^\ga\,dz=0, \ga\in\C$ (cf. (\ref{G4-1.13},b,c).
Finally, we introduce a special solution of the homogeneous
equation $lf=0$. 

A function $\varphi:(0,1)\to \R$ is called a {\it normalized solution
of $lf=0$ at $0$} (resp. $1$) if
\alpheqn
\begin{equation}
 l\varphi=0,\quad R_0 \varphi=0\quad (\mbox{\rm resp.}\;R_1\varphi=0)
\end{equation}
and
\begin{equation}
\varphi(x)= x^{\nu_0+1/2} \varphi_1(x),\quad \varphi_1(0)=1
\end{equation}
(resp. $\varphi(x)= (1-x)^{\nu_1+1/2} \varphi_1(x),\quad \varphi_1(1)=1$).
\reseteqn

It is clear that a normalized solution at $0$\quad (resp. $1$) exists and is uniquely
determined.

\label{p5}

Now we can state our main result:

\begin{theorem}\label{S4.1} Let $l$ be a regular singular Sturm--Liouville operator
as defined in \myref{G1.7}, \myref{G1.8}, \myref{G1.9}. Let $R_0, R_1$ be boundary conditions as defined
before. Then we have
$$\detz(L)=\frac{\pi\,W(\psi,\varphi)}%
    {2^{\nu_0+\nu_1} \Gamma(\nu_0+1)\Gamma(\nu_1+1)},$$
where $\varphi$ is a normalized solution of $lf=0$ at $0$ and
$\psi$ is a normalized solution of $lf=0$ at $1$.
$W(\psi,\varphi)=\psi \varphi'-\psi'\varphi$ denotes the Wronskian of
$\psi,\varphi$.
\end{theorem}

Some historical remarks are appropriate here:

For smooth potentials, S. Levit and U. Smilansky \cite{LS} showed, that
$$\detz(L)=C\; W(\psi,\varphi),$$
where $C$ is a constant depending only on the boundary condition.
This is basically the variation result Proposition \plref{S2.2} below.
T. Dreyfus and H. Dym \cite{DD} generalized this result to operators of
arbitrary order.
The first who were able to calculate the constant were Burghelea,
Friedlander and Kappeler, who calculated the determinant for
smooth operators of arbitrary order. They considered periodic
\cite{BFK1} and separated boundary conditions \cite{BFK2}.

Our method of proof is similar to \cite{BFK1,BFK2}. However, we do not use
the asymptotic expansion of $\detz(L+z)$ for large $z$, nor do we use
the theory of complex functions of a certain order. Instead, the problem
is reduced to the explicit calculation of the determinant of a single
operator. Moreover we use the well-known values $\zeta_R(0)=-1/2,
\zeta_R'(0)=-\frac 12 \log (2\pi)$ of the Riemann $\zeta$--function.

\newsection{The determinant of the regular singular model operator}

\def\sec1{2}
\def\Inu{I_\nu}
\def\Knu{K_\nu}

In this section we calculate the determinant of the Friedrichs extension
of the model operator
$l_\nu$ \myref{G1.14}, $\nu\ge 0$. Let $R_1 f= f(1)$ and put
\begin{equation}
L_\nu:= (l_\nu, R_0, R_1)=l_\nu^\cf.
\label{G3.30}
\end{equation}
The kernel, $k(x,y;z)$, of the resolvent $(L_\nu+z^2)^{-1}$ is given
in terms of the modified Bessel functions $I_\nu, K_\nu$ (cf.
\cite{BS85}) by
\begin{equation}
k(x,y;z)= (xy)^{1/2} I_\nu(xz)\Big (K_\nu (yz)-
   \frac{K_\nu(z)}{I_\nu(z)}\Inu(yz)\Big),\quad x\le y.
   \label{G3.31}
\end{equation}

\def\Lnu{L_\nu}
Moreover,  $\Lnu$ is invertible and the kernel of $\Lnu^{-1}$ is
\begin{equation}
    k_\nu(x,y):=k_\nu(x,y;0)=\casetwo{x^{\nu+1/2}(y^{-\nu+1/2}-y^{\nu+1/2})/2\nu}{\nu>0,\; x\le y,}%
      {-x^{1/2}y^{1/2}\;\log(y)}{\nu=0,\; x\le y.}
   \label{G3.1}
\end{equation}

We adopt the following notation: multiplication operators by functions
are denoted by the corresponding capital letters. For example the
multiplication operator by x is denoted by $X$.

\begin{lemma}\label{S3.1} Let $\nu>0$. Then $X^{-1}\Lnu^{-1}\in C_2(L^2[0,1])$ is
a Hilbert--Schmidt operator. Moreover $\nu\mapsto X^{-1}\Lnu^{-1}$,
$\nu>0$, is a continuous map into $C_2(L^2[0,1])$.
\end{lemma}
\begin{proof} This follows immediately from the kernel representation.
For instance, we have for $x\le y$
\begin{eqnarray*}
    2\nu x^{-1} |k_\nu(x,y)|&\le& x^{\nu-1/2}y^{-\nu+1/2}+ x^{\nu-1/2}y^{\nu+1/2}\\
       &\le& 2 x^{\nu-1/2}y^{-\nu+1/2},
\end{eqnarray*}
thus
$$\nu^2 \int_0^1\int_0^y x^{-2} |k_\nu(x,y)|^2 dx dy \le
   \int_0^1 y^{1-2\nu} \int_0^y x^{2\nu-1} dx dy = 1/4\nu<\infty.$$
The estimate for $x\ge y$ is similar and the continuity statement
is obvious.\end{proof}

\def\trlznu{\Tr((L_\nu+z^2)^{-1})}
\def\Lnuinv{(L_\nu+z^2)^{-1}}
\def\Lnumtwo{(L_\nu+z^2)^{-2}}
From this lemma we infer that $\trlznu, \nu>0$, is differentiable with
respect to $\nu$ and
\begin{equation}
\frac{d}{d\nu} \trlznu= -2\nu\; \Tr(X^{-1}\Lnumtwo X^{-1}).
\end{equation}
Namely, for $\nu_1, \nu_2$ we find
  $$(L_{\nu_1}+z^2)^{-1} - (L_{\nu_2}+z^2)^{-1}=
    (\nu_2^2-\nu_1^2)(L_{\nu_1}+z^2)^{-1}X^{-1} X^{-1} (L_{\nu_2}+z^2)^{-1}.$$
Hence by Lemma \plref{S3.1}, $\trlznu$ is differentiable as a map
from $(0,\infty)$ into the space of trace class operators and we obtain
the formula.

Now we come to the main tool for calculating the determinant of
$\Lnu$.

\begin{proposition}\label{S3.2} $\nu\mapsto T(\Lnu), \nu>0,$ is differentiable and
$$\frac{d}{d\nu} T(L_\nu)= -\log 2 - \frac{\Gamma'}{\Gamma}(\nu+1).$$
\end{proposition}
\begin{proof} First we show that we can differentiate under the integral
in \myref{G1.6}. From the preceding considerations, we conclude that
$$(0,\infty)\times [0,\infty)\ni (\nu,z)\mapsto \trlznu$$
is continuously differentiable, hence we have
\alpheqn
\begin{equation}
\frac{d}{d\nu} \int_0^1 z\; \trlznu dz = -2\nu \int_0^1 z\; \Tr(X^{-1}
  \Lnumtwo X^{-1})dz.
  \label{G2.NOa}
\end{equation}

To show that
\begin{equation} \begin{array}{cl}
&\DST\frac{d}{d\nu} \int_1^\infty z\; \trlznu-1/2 +(1/4+\nu/2) z^{-1}\, dz\\
=&\DST\int_1^\infty -2 \nu \, z\; \Tr(X^{-1}\Lnumtwo X^{-1})+1/2\, z^{-1}\,dz,
  \end{array}\label{G2.NOb}
\end{equation}
\reseteqn
it is enough to prove the estimate
\begin{equation}|-2  z\; \Tr(X^{-1}\Lnumtwo X^{-1})+\frac{1}{2\nu} z^{-1}|
  \le c |z|^{-2}
  \label{G2.2}
\end{equation}
with $c$ locally independent of $\nu$.
Then the differentiability is a consequence of the dominated convergence
theorem.

Since
$$-2 z\; (\Lnu+z^2)^{-2} = \frac{d}{dz} \Lnuinv,$$
the kernel of $-2 z\,X^{-1}(\Lnu+z^2)^{-2} X^{-1}$ is given by
$$\frac{d}{dz} x^{-1} k_\nu(x,y;z) y^{-1},$$
hence
$$-2 z\; \Tr(X^{-1}(\Lnu+z^2)^{-2} X^{-1})=
   \int_0^1 x^{-1} \frac{d}{dz} \Inu(xz)\Big(\Knu(xz)-\frac{\Knu(z)}{\Inu(z)}
     \Inu(xz)\Big)\,dx.$$
We have
\begin{eqnarray*}
     \int_0^1 x^{-1} \frac{d}{dz} \Inu(xz)\Knu(xz)dx&=&
         z^{-1} \int_0^z \frac{d}{dx} \Inu \Knu (x)dx\\
         &=&z^{-1} \Inu\Knu(z)-\frac{1}{2\nu} z^{-1}.
\end{eqnarray*}
Furthermore
\begin{eqnarray*}
  && \int_0^1 x^{-1} \frac{d}{dz} \Inu(xz)^2 \frac{\Knu(z)}{\Inu(z)} dx\\
  &=&  \frac{\Knu(z)}{\Inu(z)} \int_0^1 x^{-1} \frac{d}{dz} \Inu(xz)^2 dx
     + \int_0^1 x^{-1} \Inu(xz)^2 dx \frac{d}{dz} \frac{\Knu(z)}{\Inu(z)}\\
   &=& z^{-1} \Inu\Knu(z) + \int_0^z x^{-1} \Inu(x)^2 dx
       \frac{d}{dz} \frac{\Knu(z)}{\Inu(z)},
\end{eqnarray*}
hence we find
\begin{eqnarray*}
   -2 z\; \Tr(X^{-1}\Lnumtwo X^{-1})+\frac{1}{2\nu} z^{-1}&=&
     -\int_0^z x^{-1} \Inu(x)^2 dx \frac{d}{dz} \frac{\Knu(z)}{\Inu(z)}.
\end{eqnarray*}

Now the estimate \myref{G2.2} follows from the asymptotics of the
modified Bessel functions \cite[7.23]{Watson}
\alpheqn
\begin{eqnarray}
    \Inu(x)&=&\DST\frac{1}{\sqrt{2\pi x}} e^x(1+ O(x^{-2})),\; x\to \infty,
        \label{G3.32a}\\
\DST    \Knu(x)&=&\DST\sqrt{\frac{\pi}{2x}} e^{-x}(1+ O(x^{-2})),\; x\to \infty,
   \label{G3.32b}
\end{eqnarray}
\reseteqn
which can be differentiated with respect to
$x$ and are locally uniform in $\nu$.

In view of Theorem \plref{S1.1}, \myref{G1.38}, (\plref{G2.NOa},b) we have
proved that $T(\Lnu)$ is differentiable and
\begin{eqnarray}
   \frac{d}{d\nu} T(\Lnu)&=& 4\nu \int_0^1 z\; \Tr(X^{-1}
     (\Lnu+z^2)^{-2} X^{-1}) dz\nonumber \\
     &&+ 2 \int_1^\infty 2 \nu z\; \Tr(X^{-1}(\Lnu+z^2)^{-2} X^{-1})
       -1/2\, z^{-1}\, dz\nonumber \\
     &=& 4\nu \regint_0^\infty z\; \Tr(X^{-1}(\Lnu+z^2)^{-2} X^{-1}) dz.
\end{eqnarray}

Furthermore, using \myref{G3.31}, we find
\begin{eqnarray*}
 &&  4\nu \regint_0^\infty z\; \Tr(X^{-1}(\Lnu+z^2)^{-2} X^{-1}) dz\\
 &=& -2\nu \regint_0^\infty \int_0^1 x^{-1} \frac{d}{dz} \Inu(xz)\Big(
    \Knu(xz)- \frac{\Knu(z)}{\Inu(z)}\Inu(xz)\Big) dx\\
 &=& -2\nu\regint_0^\infty z^{-1} \Inu\Knu(z)-\frac{1}{2\nu} z^{-1} -
   \frac{d}{dz}\left( \frac{\Knu(z)}{\Inu(z)}\int_0^1 x^{-1} \Inu(xz)^2 dx\right) dz\\
 &=& -2\nu\regint_0^\infty z^{-1} \Inu\Knu(z)dz + 2\nu
   \regint_0^\infty \frac{d}{dz}\left( \frac{\Knu(z)}{\Inu(z)}\int_0^z
     x^{-1} \Inu(x)^2 dx\right) dz\\
 &=:& -2\nu I_1 +2\nu I_2.
\end{eqnarray*}
The first integral is well--known (cf. e.g. \cite[p. 418]{BS87}).
One has, more generally,
\begin{equation}
\regint_0^\infty x^s \Inu\Knu(x) dx = \frac{\Gamma(\frac{s+1}{2})\Gamma(-s/2)\Gamma(\nu+\frac{s+1}{2})}%
   {4\sqrt{\pi} \Gamma(\nu-\frac{s-1}{2})}.
\end{equation}
Since the right hand side has a simple pole at $s=-1$, we find using
the Mellin--transform definition of $\reginttext$
\begin{eqnarray*}
I_1 &=& \frac{d}{ds}\big|_{s=-1} (s+1)
   \frac{\Gamma(\frac{s+1}{2})\Gamma(-s/2)\Gamma(\nu+\frac{s+1}{2})}%
   {4\sqrt{\pi} \Gamma(\nu-\frac{s-1}{2})}\\
   &=&\frac{\log 2}{2\nu} +\frac{1}{4\nu} \Big(\frac{\Gamma'}{\Gamma}(\nu)
      +\frac{\Gamma'}{\Gamma}(\nu+1)\Big).
\end{eqnarray*}
Since
$$\frac{\Knu}{\Inu}(z) \int_0^z x^{-1} \Inu(x)^2dx =O(z^{-2}),\quad z\to \infty,$$
$I_2$ is actually a regular integral and we find
$$I_2= -\lim_{z\to 0}\frac{\Knu}{\Inu}(z) \int_0^z x^{-1} \Inu(x)^2dx
   =-(1/2\nu)^2.$$
Thus we end up with
\begin{eqnarray*}
    \frac{d}{d\nu} T(L_\nu)&=& -\log 2 -\frac 12
       \Big(\frac{\Gamma'}{\Gamma}(\nu)
    +\frac{\Gamma'}{\Gamma}(\nu+1)\Big)-\frac{1}{2\nu}\\
    &=&-\log 2-\frac{\Gamma'}{\Gamma}(\nu+1).
\end{eqnarray*}
\end{proof}

An immediate consequence is the

\begin{theorem}\label{S3.3} Let $L_\nu, \nu\ge 0,$ be the operator defined
in \myref{G3.30}. Then the $\zeta$--regularized determinant of
$L_\nu$ is given by
$$\detz L_\nu=\frac{\sqrt{2\pi}}{2^\nu \Gamma(\nu+1)}.$$
\end{theorem}
\begin{proof} From the preceding proposition we infer for $\nu>0$
$$\frac{d}{d\nu} \log \detz L_\nu=-\frac{d}{d\nu}(\nu\, \log 2 +\log \Gamma
  (\nu+1))=\frac{d}{d\nu} \log \frac{\sqrt{2\pi}}{2^\nu \Gamma(\nu+1)}, $$
hence it suffices to check the formula for $\nu=1/2$ and $\nu=0$. 
$L_{1/2}$ is just $-\frac{d}{d x^2}$ with Dirichlet boundary conditions,
$$\spec(L_{1/2})=\{ \,n^2 \pi^2\,|\, n\in\N\,\},$$
thus
$$\zeta_{L_{1/2}}(s)= \pi^{-2s} \zeta_R(2s),$$
where $\zeta_R$ denotes the Riemann zeta--function. In view of the well--known
formulas
$\zeta_R(0)=-1/2,\quad \zeta'_R(0)=-\frac 12 \log(2\pi)$ we find
\begin{eqnarray*}
    \log\detz(L_{1/2})&=&-\zeta_{L_{1/2}}'(0)=2\,\log \pi\, \zeta_R(0)-2\zeta_R'(0)\\
       &=& \log 2 =\log\frac{\sqrt{2\pi}}{\sqrt{2} \Gamma(3/2)}.
\end{eqnarray*}
To prove the result for $\nu=0$ it is enough to show that
$\nu\mapsto \detz(L_\nu)$ is continuous at $\nu=0$. Similar to
the argument \myref{G2.2} in the proof of Proposition \plref{S3.2},
it suffices to prove the estimate
\begin{equation}
|z\, \Tr((L_\nu+z^2)^{-1})-1/2 +\frac 12 (\nu+1/2) z^{-1}|\le c\, |z|^{-2}
  \label{G3.33}
\end{equation}
with $c$ locally independent of $\nu$. Then continuity is a consequence
of the dominated convergence theorem. The estimate \myref{G3.33}
is a consequence of Theorem \plref{S1.1}. That the constant $c$ is indeed
locally independent of $\nu$
follows easily from the asymptotic relations (\plref{G3.32a},b).\end{proof}

\newsection{Variation formulas}

\begin{lemma}\label{S4-3.0} Let $L_\nu, \nu\ge 0,$ be
the operator defined in \myref{G3.30}. Then for
$\delta>0, z\ge0$ the operator $X^{\delta-1}(L_\nu+z^2)^{-1/2}$
is Hilbert--Schmidt and we have the estimate
\begin{eqnarray*}
      \|X^{\delta-1}(L+z^2)^{-1/2}\|_{C_2}=
        \casethree{O(z^{-\delta})}{0<\delta<1/2,}{O(z^{-1/2}\log z)}{\delta=1/2,}{%
           O(z^{-1/2})}{\delta>1/2,}
\end{eqnarray*}
as $z\to\infty$.
\end{lemma}
\begin{proof} Introducing
the first order operator
\begin{equation}
  D_\nu:=\frac{d}{dx}+(\nu-1/2)X^{-1}
  \label{G4-3.0}
\end{equation}
one checks that $L_\nu=D_{\nu,\max}^tD_{\nu,\min}$
(cf. \cite[Lemma 3.1]{BL93}). Moreover, since $L_\nu$ is the Friedrichs
extension of $l$, the domain of $L_\nu^{1/2}$ is the completion
of $\cinfz{0,1}$ with respect to the norm
\begin{eqnarray*}
     \|f\|_{L_\nu^{1/2}}^2&=& (f,f)+(L_\nu^{1/2} f,L_\nu^{1/2} f)=
        (f,f)+(L_\nu f,f)\\
        &=&(f,f)+(D_\nu^t D_\nu f,f)=(f,f)+(D_\nu f,D_\nu f).
\end{eqnarray*}
But since the latter is the square of the graph norm of $D_\nu$
we find $\cd(L_\nu^{1/2})=\cd(D_{\nu,\min})$.

From \cite[Lemma 2.1]{BS88} we have
for $f\in\cd(D_{\nu,\min})=\cd(L_\nu^{1/2})$
\begin{equation}
    |f(x)|\le c |x\, \log x|^{1/2} \|D_{\nu,\min} f\|=c|x\,\log x|^{1/2}
    \|L_\nu^{1/2}f\|.
    \label{G2.40}
\end{equation}
Now let $\mltilde k(x,y;z)$ be the kernel of $(L+z^2)^{-1/2}$. Then
\myref{G2.40} implies in view of the Theorem of Riesz
$$\int_0^1|\mltilde k(x,y;z)|^2 dy \le c\, x\, |\log x|$$
with $c$ independent of $z\ge 0$. Thus $X^{\delta-1}(L_\nu+z^2)^{-1/2}$ is
a Hilbert--Schmidt operator with Hilbert--Schmidt norm uniformly bounded
for $z\ge 0$.

This proves that $X^{\delta-1}(L_\nu+z^2)^{-1}X^{\delta-1}$ is trace
class for $z\ge 0$. Moreover, since $X^{\delta-1}(L_\nu+z^2)^{-1}X^{\delta-1}\ge 0$
we find
\begin{eqnarray*}
  && \|X^{\delta-1}(L_\nu+z^2)^{-1/2}\|_{C_2}^2=
        \Tr (X^{\delta-1}(L_\nu+z^2)^{-1}X^{\delta-1})\\
  &&\quad =\int_0^1 x^{2\delta-1} I_\nu K_\nu(xz)dx - \frac{K_\nu}{I_\nu}(z)
           \int_0^1 x^{2\delta-1} I_\nu(xz)^2 dx\\
  &&\quad = z^{-2\delta}\left(\int_0^z u^{2\delta-1} I_\nu K_\nu (u)du
           - \frac{K_\nu}{I_\nu}(z)\int_0^z u^{2\delta-1}I_\nu(u)^2du
           \right).
\end{eqnarray*}
Now the assertion follows easily from the well--known asymptotics
of $I_\nu(x), K_\nu(x)$ as $x\to 0$
\begin{equation}
   I_\nu(x)\sim c x^\nu,\quad K_\nu(x)\sim\casetwo{c\, x^{-\nu}}{\nu>0,}{c\, \log x}{\nu=0,}
   \label{G2.41}
\end{equation}
and the asymptotics (\plref{G3.32a},b) as $x\to\infty$.
\end{proof}

\label{page14}

Next we sketch the construction of the resolvent of general
$L=(l,R_0,R_1)$ (cf. \cite[Sec. 4]{BS85}). Let
\begin{equation}
l=-\frac{d^2}{d x^2} + q(x)
\end{equation}
be a regular singular Sturm--Liouville operator as defined in
\myref{G1.7}, \myref{G1.8}, \myref{G1.9}. Let $L=(l,R_0,R_1)$
be a self--adjoint
extension. We consider the case that $q$ is not continuous at both
ends. The other cases are easier. We choose cut--off functions
$\varphi, \mltilde\varphi, \psi, \mltilde\psi\in\cinf{[0,1]}$ as
follows:
\begin{equation}\begin{array}{ll}
    \supp\varphi\subset [0,1/3], & \supp\psi\subset [2/3,1],\\
    \varphi|[0,1/6]=1,& \psi|[5/6,1]=1,\\
    \mltilde\varphi|[0,1/3]=1,& \mltilde\psi|[2/3,1]=1,\\
    \mltilde\varphi+\mltilde\psi=1.
\end{array}\label{G4-3.1}
\end{equation}
Then we can write
\begin{equation}
  l=: -\frac{d^2}{dx^2}+\varphi (\nu_0^2-1/4)X^{-2} +\psi (\nu_1^2-1/4)
      (1-X)^{-2}+\mltilde q(x)\label{G4-3.2}
\end{equation}
and
\begin{equation}
   \mltilde q(x)=O(x^{1/N-2}), x\to 0,\quad \mltilde q(x)=O((1-x)^{1/N-2}), x\to 1.
   \label{G4-3.3}
\end{equation}

For $\nu\ge 0$ we denote by $\mltilde L_\nu$ the Friedrichs extension
of the operator
$$-\frac{d^2}{dx^2}+(1-x)^{-2}(\nu^2-1/4)$$
on $\cinfz{0,1}\subset L^2[0,1]$ (cf. \myref{G3.30}).

Now we put
\begin{equation}
  G(z):= \mltilde\varphi (L_{\nu_0}+z)^{-1}+\mltilde\psi
          (\mltilde L_{\nu_1}+z)^{-1}.
          \label{G4-3.4}
\end{equation}
In view of \myref{G1.19}, \myref{G1.22}, $G(z)$ maps $L^2[0,1]$
into the domain of $L=(l,R_0,R_1)$. Furthermore,
\begin{eqnarray}
     (L+z)G(z)&=& I+ [L,\mltilde\varphi](L_{\nu_0}+z)^{-1} +
                                  [L,\mltilde\psi](\mltilde L_{\nu_1}+z)^{-1} \nonumber\\
              &&+\mltilde\varphi(\mltilde q+(\varphi-1)(\nu_0^2-1/4)X^{-2})(L_{\nu_0}+z)^{-1}\nonumber\\
              &&+\mltilde\psi(\mltilde q+(\psi-1)(\nu_1^2-1/4)(1-X)^{-2})(\mltilde L_{\nu_1}+z)^{-1}
                \label{G4-3.5}\\
             &=:&I+R(z),\nonumber
\end{eqnarray}
where $[L,\mltilde\varphi]=-\mltilde\varphi''-2\mltilde\varphi'\frac{d}{dx}$.

\begin{lemma}\label{S4-3.1} We use the notation introduced before.

{\rm 1.} For any function $\chi\in\cinfz{0,1}$ we have for $|z|\ge z_0$
$$\|\chi \frac{d}{dx} (\stackrel{(\sim)}{L_\nu}+z)^{-1/2}\|_{C_2}\le C(\chi,\nu,z_0).$$

{\rm 2.} For $|z|\ge z_0$ we have $\|R(z)\|\le C(z_0) |z|^{-1/2}$,
where the constant $C(z_0)$ depends only on $\varphi, \mltilde\varphi, \psi,
\mltilde \psi$ and
$$\DST \sup_{0<x<1}\TST X^{2-1/N}(1-X)^{2-1/N}\mltilde q(x).$$
Furthermore, for $|z|$ large,
$$(L+z)^{-1}=G(z) \sum_{n=0}^\infty (-1)^n R(z)^n.$$
\end{lemma}
\begin{proof} 1. Since $l$ is elliptic of order $2$, we have $H^1_{\rm loc}(0,1)
\subset \cd(L_\nu^{1/2})$ from which we reach
the conclusion immediately.

2. In view of \myref{G4-3.5} we only have to prove the estimate
$\|R(z)\|\le C(z_0) |z|^{-1/2}$. But this is an easy consequence of
Lemma \plref{S4-3.0}, \myref{G4-3.5} and the proven first part of this
lemma.\end{proof}

\begin{lemma}\label{S4-3.2} We use the notation of page {\rm \pageref{page14}.}
For $\delta>0, z>0$ the operator
$X^{\delta-1}(L+z^2)^{-1/2}$ is Hilbert--Schmidt and
\begin{eqnarray*}
      \|X^{\delta-1}(L+z^2)^{-1/2}\|_{C_2}=
        \casethree{O(z^{-\delta})}{0<\delta<1/2,}{O(z^{-1/2}\log z)}{\delta=1/2,}{%
           O(z^{-1/2})}{\delta>1/2,}
\end{eqnarray*}
as $z\to \infty$. Again, the O--constant depends only on
$\varphi, \mltilde\varphi, \psi, \mltilde\psi$ and
$$\DST \sup_{0<x<1}\TST X^{2-1/N}(1-X)^{2-1/N}\mltilde q(x).$$
The same estimate holds for $\|(1-X)^{\delta-1}(L+z^2)^{-1/2}\|_{C_2}$.

If $L$ is invertible, then $X^{\delta-1} L^{-1}$ is Hilbert--Schmidt, too.
\end{lemma}
\begin{proof} Lemma \plref{S4-3.0}, Lemma \plref{S4-3.1} and the formula
\begin{eqnarray}
   (L+z^2)^{-1/2}&=& \frac{1}{\pi}\int_0^\infty \gl^{-1/2} (L+z^2+\gl)^{-1} d\gl
                    \nonumber\\
          &=& \frac{1}{\pi} \int_0^\infty \gl^{-1/2} G(z^2+\gl)
              \sum_{n=0}^\infty (-1)^n R(z^2+\gl)^n d\gl,\label{G4-3.6}
\end{eqnarray}
which holds for $z>0$ large enough, imply the assertion for $z\ge z_0$.
If $z\in \C$ with $L+z^2$ invertible, then we conclude from
$$X^{\delta-1}(L+z^2)^{-1/2}=X^{\delta-1}(L+z_0^2)^{-1/2}\big[(L+z_0^2)^{1/2}(L+z^2)^{-1/2}\big]$$
that the operator $X^{\delta-1}(L+z^2)^{-1/2}$ is Hilbert--Schmidt,
too.\end{proof}

Now we introduce smooth families of operators.
$$l_t=-\frac{d^2}{dx^2}+q_t(x)$$
is said to be a smooth family of operators if
$$q_t(x)=a_0(t,x^{1/N}) x^{-2}=(1-x)^{-2} a_1(t,(1-x)^{1/N})$$
with
smooth functions $a_0\in \cinf{I\times[0,1)}, a_1\in\cinf{I\times(0,1]}$,
$I$ some interval.

\begin{proposition}\label{S2.2} Let $l_t$ be a smooth family of operators with
$\nu_0, \nu_1$ independent of $t$. Let $R_0, R_1$ be fixed boundary
conditions independent of $t$. Moreover let $\varphi_t, \psi_t$ be normalized
solutions of $l_t f=0$ at $0$ resp. $1$.
If $L_t=(l_t, R_0, R_1)$ is invertible, then
$T(L_t)$ is smooth and we have the variation formula
$$\frac{d}{dt} T(L_t)=\frac{d}{dt} \log W(\psi_t,\varphi_t).$$
Here, $W(\psi_t,\varphi_t)=\psi_t\varphi_t'-\psi_t'\varphi_t$ denotes
the Wronskian of $\psi_t, \varphi_t$.
\end{proposition}

\begin{remark} This Proposition is essentially the result of \cite{LS} and our
proof is an adaption of their proof to our more general setting.
\end{remark}

\begin{proof} Since $\nu_0, \nu_1$ are independent of $t$, we have the
estimate
\begin{equation}
|\pl_t q_t(x)|\le c\, x^{1/N-2}(1-x)^{1/N-2}
  \label{G4-3.7}
\end{equation}
with $c$ locally independent of $t$.

We would like to apply the formula
$\frac {d}{dt} T(L_t)= \Tr((\pl_t q_t)  L_t^{-1}).$ However, as the
referee pointed out to the author, the operator $(\pl_t q_t)  L_t^{-1}$
need not be of trace class. But, in view of Lemma \plref{S4-3.2}
the operator
$$X^{1-1/2N}(1-X)^{1-1/2N}(\pl_t q_t)  L_t^{-1}X^{1/2N-1}(1-X)^{1/2N-1}$$
is trace class and the kernels of this operator and
$(\pl_t q_t) L_t^{-1}$ 
coincide on the diagonal.
We introduce the abbreviation $\omega(x):=x(1-x)$ and denote by $\Omega$
the operator of multiplication by $\omega$.

To make the preceding consideration rigorous we recall from \myref{G1.38}
(note that $L_t$ is assumed to be invertible)
\begin{eqnarray*}
     T(L_t)&=& -2 \regint_0^\infty z\, \Tr(L_t+z^2)^{-1} dz\\
        &=& -2 \regint_0^1 z\, \Tr(L_t+z^2)^{-1} dz\\
        &&-2 \regint_1^\infty z\,[\Tr(L_t+z^2)^{-1} -a_0 z^{-1}
           -c_0 z^{-2}]dz.
\end{eqnarray*}
$a_0, c_0$ are independent of $t$ in view of Theorem \plref{S1.1}.
Formal differentiation under the integral gives
\begin{equation}
\frac{d}{dt} T(L_t)= 2 \int_0^\infty z\,\Tr ((L_t+z^2)^{-1}
  (\pl_t q_t) (L_t+z^2)^{-1})dz
  \label{G4-3.8}
\end{equation}
To justify this formula we estimate the integrand using Lemma \plref{S4-3.2}
\begin{eqnarray*}
  &&  |\Tr((L_t+z^2)^{-1}(\pl_t q_t)(L_t+z^2)^{-1})|\\
  &\le& \|(L_t+z^2)^{-1}(\pl_t q_t)\Omega^{1-1/2N}\|_{C_2}
     \|\Omega^{1/2N-1}(L_t+z^2)^{-1}\|_{C_2}\\
  &\le& \|(L_t+z^2)^{-1/2}\|^2 \|(L_t+z^2)^{-1/2}(\pl_t q_t)\Omega^{1-1/2N}\|_{C_2}
     \|\Omega^{1/2N-1}(L_t+z^2)^{-1/2}\|_{C_2}\\
 &=&O(z^{-2-1/N})
\end{eqnarray*}
where the O--constant is locally independent of $t$.
Now \myref{G4-3.8} follows from the dominated convergence theorem.
We continue starting from \myref{G4-3.8} and find
\begin{eqnarray*}
    \frac{d}{dt} T(L_t)&=& 2 \int_0^\infty z\, \Tr((L_t+z^2)^{-1}
       \Omega^{1/2N-1} \Omega^{1-1/2N}(\pl_t q_t)(L_t+z^2)^{-1}) dz\\
    &=&2\int_0^\infty z\, \Tr(\Omega^{1-1/2N}(\pl_t q_t)(L_t+z^2)^{-2}
        \Omega^{1/2N-1})dz\\
    &=& -\int_0^\infty \frac{d}{dz} \Tr(\Omega^{1-1/2N}(\pl_t q_t)
               (L_t+z^2)^{-1}\Omega^{1/2N-1})dz\\
   &=&\Tr(\Omega^{1-1/2N}(\pl_t q_t)L_t^{-1}\Omega^{1/2N-1})\\
   &=& \int_0^1 ((\pl_t q_t)L_t^{-1})(x,x) dx,
\end{eqnarray*}
which morally is $\Tr((\pl_t q_t)L_t^{-1})$ although the latter in general
does not exist.

The kernel of $L_t^{-1}$ is given by
\begin{equation}
k(x,y)=W(\psi,\varphi)^{-1}\varphi(x) \psi(y),\quad x\le y.
  \label{G2.42}
\end{equation}
Note that $W(\psi,\varphi)\not=0$ since $L_t$ is assumed to be invertible.
Differentiating the formula $\varphi''=q_t \varphi$ with respect to $t$ gives
$$\pl_t\varphi''=(\pl_t q_t) \varphi+q_t \pl_t \varphi$$
and hence
\begin{eqnarray*}
    (\pl_t q_t)\varphi\psi&=&(\pl_t\varphi)''\psi -q_t (\pl_t \varphi) \psi=
         (\pl_t \varphi)''\psi-(\pl_t \varphi)\psi''\\
         &=&\frac{d}{dx}\big[(\pl_t \varphi)'\psi -(\pl_t \varphi) \psi'\big]\\
         &=&\frac{d}{dx} W(\psi,\pl_t\varphi).
\end{eqnarray*}
Thus we find
\begin{eqnarray*}
    \frac{d}{dt}T(L_t)&=&\int_0^1 ((\pl_t q_t) L_t^{-1})(x,x) dx\\
       &=&W(\psi,\varphi)^{-1}\int_0^1\frac{d}{dx}W(\psi,\pl_t\varphi)(x) dx\\
       &=&W(\psi,\varphi)^{-1}\big[ W(\psi,\pl_t\varphi)(1)-
          W(\psi,\pl_t\varphi)(0)\big].
\end{eqnarray*}

Since $\varphi$ is normalized at 0 and $\nu_0$ is constant, we have
$$\pl_t \varphi(x)= O(x^{\nu_0+1/2+1/N}), \quad x\to 0$$
and
$$\pl_t \varphi'(x)=O(x^{\nu_0-1/2+1/N}), \quad x\to 0,$$
which implies immediately
 $$W(\psi,\pl_t\varphi)(x)=O(x^{1/N}\,\log(x)),\quad x\to 0$$
and thus $W(\psi,\pl_t\varphi)(0)=0$.

Reversing the roles of $\varphi,\psi$ we find $W(\pl_t \psi,\varphi)(1)=0$.
Summing up we have
\begin{eqnarray*}
       \frac{d}{dt} T(L_t)&=& W(\psi,\varphi)^{-1}\big[W(\psi,\pl_t\varphi)(1)
          +W(\pl_t\psi,\varphi)(1)\big]\\
          &=&W(\psi,\varphi)^{-1}\frac{d}{dt} W(\psi,\varphi)=
             \frac{d}{dt} \log \,W(\psi,\varphi).
\end{eqnarray*}
\end{proof}

The next Proposition is basically \cite[Proposition 3.2]{BFK2}. The fact
that $q$ may be singular at $0$ causes no essential new difficulty.
To make the exposition self--contained we include a proof.

\begin{proposition}\label{S2.3} Assume that $q$ is continuous at $1$ and let
$R_{1,t}f= f'(1)+a(t) f(1)$ be a smooth family of boundary operators of
order $1$.
Assume that $L_t=(l,R_0,R_{1,t})$ is invertible. Then we have
$$\frac{d}{dt} T(L_t)= \frac{d}{dt} \log\, W(\psi_t,\varphi_t),$$
where $\varphi_t,\psi_t$ are as in Proposition \plref{S2.2}.
\end{proposition}
\begin{proof} For simplicity, throughout this proof we are going to write
$\varphi, \psi$ instead of $\varphi_t, \psi_t$.

Let $g\in L^2[0,1]$ and consider $f_t:=L_t^{-1} g$. Differentiation
of $l f_t=g$ and $R_{1,t}f_t=0$ gives
$$l \pl_t f_t=0,\quad R_{1,t} \pl_t f=-a'(t) f_t(1).$$
Since $L_t$ is invertible, we have $R_{1,t} \varphi\not=0$.
Now note that $\varphi$ is independent of $t$ and
\begin{equation}
  W(\psi,\varphi)(1)=\psi(1) \varphi'(1)-\psi'(1)\varphi(1)=
    \varphi(1) a(t)+\varphi'(1)=R_{1,t}\varphi.
    \label{G2.43}
\end{equation}

Now consider
$$\pl_t f_t+a'(t) W(\psi,\varphi)^{-1} f_t(1)\varphi=u.$$
We find
$$lu=0,\quad R_0u=R_1u=0$$
and again since $L_t$ is invertible,
$$\pl_t f_t=-a'(t)W(\psi,\varphi)^{-1} f_t(1) \varphi.$$
Thus $\pl_t L_t^{-1}$ is actually a rank one operator (see \myref{G2.42}):
$$(\pl_t L_t^{-1}g)(x)=-a'(t) W(\psi,\varphi)^{-2} \varphi_t(x)
   \int_0^1 \varphi_t(y) g(y) dy.$$
Now let $\varphi_t(x,z), \psi_t(x,z)$ be the corresponding solutions
for $L_t+z^2$. Then we find
\begin{eqnarray*}
    \pl_t \Tr((L_t+z^2)^{-1})&=& -a'(t) W(\psi,\varphi)^{-2}\int_0^1 \varphi_t(y,z)^2dy\\
    &=&-a'(t) \int_0^1 (L_t+z^2)^{-1}(1,y)^2 dy\\
    &=&-a'(t) (L_t+z^2)^{-2}(1,1)\\
    &=& a'(t) \frac{1}{2z} \frac{d}{dz} (L_t+z^2)^{-1}(1,1)
\end{eqnarray*}
and we reach the conclusion using \myref{G2.43} and \myref{G1.38}
\begin{eqnarray*}
   \frac{d}{dt} T(L_t)&=&-a'(t)\regint_0^\infty \frac{d}{dz} (L_t+z^2)^{-1}(1,1)dz\\
   &=& a'(t) (L_t+z^2)^{-1}(1,1)\\
   &=&a'(t) W(\psi,\varphi)^{-1} \varphi(1)
\end{eqnarray*}
and we are done.\end{proof}

\begin{proposition}\label{S2.4} Let $q_\nu(x):=(\nu^2-1/4)/x^2+q_1(x), \nu>0,$ with
$\supp(q_1)\subset [\eps,1]$ for some $\eps>0$ and put
$\mltilde l_\nu:=-\frac{d^2}{dx^2}+q_\nu$. Let $R_0$ be as in \myref{G1.19}
and choose $R_1$ fixed. Moreover let $\varphi_\nu, \psi_\nu$ be normalized
solutions of $\mltilde l_\nu f=0$ at $0$ resp. $1$. If $\mltilde L_\nu:=(\mltilde l_\nu,R_0,R_1)$
is invertible, then $T(\mltilde L_\nu)$ is smooth and we have the
variation formula
$$\frac{d}{d\nu} T(\mltilde L_\nu)= \frac{d}{d\nu} \log\, W(\psi_\nu,\varphi_\nu)+
       \frac{d}{d\nu} T(L_\nu)=\frac{d}{d\nu} \log\,
       \frac{W(\psi_\nu,\varphi_\nu)}{2^\nu \Gamma(\nu+1)}.$$
\end{proposition}
\begin{remark} Note that $q_1(x)$ may be singular at $1$ (cf. \myref{G1.8}).
This is the reason why this proposition is needed. If $q_1$ is smooth
on $[0,1]$ we can just apply Proposition \plref{S2.2} and deform
$q_1$ to $0$.
\end{remark}
       
\begin{proof} The resolvent expansion Lemma \plref{S4-3.2} shows that
the estimate \myref{G2.2} holds for $\mltilde L_\nu$, too. Then as in the first
part of the proof of Proposition \plref{S3.2} one infers that
$T(\mltilde L_\nu)$ is smooth and
\begin{equation}
\frac{d}{d\nu} T(\mltilde L_\nu)=4\nu\,\regint_0^\infty z\, \Tr(X^{-1}
   (\mltilde L_\nu+z^2)^{-2} X^{-1}) dz.
   \label{G4-3.10}
\end{equation}
Now, since $\supp q_1\subset [\eps,1]$ we have
$$x^{-2} (\mltilde L_\nu+z^2)^{-1}(x,x)- x^{-2} k_\nu(x,x;z)\in L^1[0,1].$$
To see this we consider $\varphi_\nu, \psi_\nu$ defined in
(\ref{G4-1.25},b). We have
\begin{eqnarray*}
    k_\nu(x,y)&=& \varphi_\nu(x) \psi_\nu(y),\quad x\le y\\
      &=& x^{\nu+1/2}(y^{-\nu+1/2}-y^{\nu+1/2}).
\end{eqnarray*}
Since $q_1|[0,\eps]=0, \varphi_\nu, \psi_\nu|[0,\eps]$ is also a fundamental
system of solutions of the homogeneous equation $\mltilde L_\nu f=0$ in the
interval $[0,\eps]$. Let $\mltilde\psi$ be the unique function with
$$ \mltilde L_\nu \mltilde\psi=0, \quad R_1 \mltilde\psi=0, W(\mltilde\psi|[0,\eps],\varphi_\nu|[0,\eps])=1.$$
Then for $0\le x\le y\le \eps$ the kernel of $\mltilde L_\nu^{-1}$ is given by
$$ \mltilde L_\nu^{-1}(x,y)=\varphi_\nu(x) \mltilde\psi(y).$$
Moreover, since $(L_\nu\mltilde\psi)(x)=0, x\le \eps$ there exist
constants $a,b$ such that
$$\mltilde\psi|[0,\eps]= a\varphi_\nu|[0,\eps]+b \psi_\nu|[0,\eps].$$
Furthermore
$$1=W(\mltilde\psi|[0,\eps],\varphi_\nu|[0,\eps])=b W(\psi_\nu,\varphi_\nu)=b$$
and hence we find for $x\le\eps$
\begin{eqnarray*}
     (\mltilde L_\nu^{-1}-k_\nu)(x,x)&=& \varphi_\nu(x)( a \varphi_\nu(x)+\psi_\nu(x))
         -\varphi_\nu(x) \psi_\nu(x)\\
         &=&(a-1) \varphi_\nu(x)^2=O(x^{2\nu+1}),\quad x\to 0
\end{eqnarray*}
which shows \myref{G4-3.10} for $z=0$. For arbitrary $z$ the proof is similar.

Thus it makes sense to abbreviate
\begin{eqnarray*}
&&\int_0^1 x^{-2} \left[(\mltilde L_\nu+z^2)^{-1}(x,x)- k_\nu(x,x;z)\right]
   dx\\
   &=:& \Tr(X^{-1}\left[(\mltilde L_\nu+z^2)^{-1}-(L_\nu+z^2)^{-1}\right]
    X^{-1}),
\end{eqnarray*}
although we do not claim that $X^{-1}\left[(\mltilde L_\nu+z^2)^{-1}-(L_\nu+z^2)^{-1}\right]
    X^{-1}$ is really trace class.

Now we have
\begin{eqnarray*}
     \frac{d}{d\nu} T(\mltilde L_\nu)&=&-2\nu \regint_0^\infty \frac{d}{dz}
      \Tr(X^{-1}\left[(\mltilde L_\nu+z^2)^{-1}-(L_\nu+z^2)^{-1}\right]
    X^{-1})dz + \frac{d}{d\nu} T(L_\nu)\\
    &=&2\nu\, \Tr(X^{-1}(\mltilde L_\nu^{-1}-L_\nu^{-1})X^{-1})+\frac{d}{d\nu}
      T(L_\nu)\\
    &=&2\nu \regint_0^1 x^{-2} \mltilde L_\nu^{-1}(x,x) dx - 2\nu \regint_0^1
       x^{-2} L_\nu^{-1}(x,x) dx +\frac{d}{d\nu}T(L_\nu).
\end{eqnarray*}
Using \myref{G3.1} we find
$$ -2\nu \regint_0^1 x^{-2} L_\nu^{-1}(x,x)dx= - \regint_0^1 x^{-1}-x^{2\nu-1} dx =\frac{1}{2\nu},$$
and as in the proof of Proposition \plref{S2.2} we infer
\begin{eqnarray*}
2\nu \regint_0^1 x^{-2} \mltilde L_\nu^{-1}(x,x)dx&=& W(\psi_\nu,\varphi_\nu)^{-1}
  \left[W(\psi_\nu,\frac{d}{d\nu}\varphi_\nu)(1)-\LIM_{x\to 0}
     W(\psi_\nu,\frac{d}{d\nu}\varphi_\nu)(x)\right]\\
     &=&W(\psi_\nu,\varphi_\nu)^{-1}
  \left[\frac{d}{d\nu}W(\psi_\nu,\varphi_\nu)-\LIM_{x\to 0}
     W(\psi_\nu,\frac{d}{d\nu}\varphi_\nu)(x)\right].
\end{eqnarray*}
A direct calculation shows that
$$W(\psi_\nu,\varphi_\nu)^{-1} W(\psi_\nu,\frac{d}{d\nu}\varphi_\nu)(x)
  = \log\,x +\frac{1}{2\nu}+O(x^{2\nu}),\quad x\to 0,$$
hence
$$\LIM_{x\to 0}
W(\psi_\nu,\varphi_\nu)^{-1} W(\psi_\nu,\frac{d}{d\nu}\varphi_\nu)(x)
=\frac{1}{2\nu}$$
and the result is proved.\end{proof}

\newsection{Proof of the main result and examples}

\noindent
{\bf Proof of Theorem \ref{S4.1}}\quad\enspace
If $L$ is not invertible, then $\varphi$ satisfies both boundary
conditions, hence $W(\psi,\varphi)=0$. So we may assume that $L$ is
invertible.
For $z\in\C$ consider $L+z$ and let $\varphi(\cdot,z), \psi(\cdot,z)$
be the corresponding normalized solutions. Then $\detz(L+z)$ and
$W(\psi(\cdot,z),\varphi(\cdot,z))$ are holomorphic functions in $\C$
and in view of Proposition \plref{S2.2} these functions have the same
logarithmic derivative. Hence it suffices to prove the formula for $L+z$
and $\Re z$ large.

We can deform the potential $q(x)$, such that
$$ q(x)=\casetwo{(\nu_0^2-1/4) x^{-2}}{x\le \eps,}{(\nu_1^2-1/4)(1-x)^{-2}}{x\ge 1-\eps,}$$
and again by Proposition \plref{S2.2} it suffices to prove the result
for these potentials and $\Re z$ large. If $\Re z$ is large
enough we apply Proposition \plref{S2.4} and deform $\nu_0$ and
$\nu_1$ to $\pm 1/2$ leaving a potential $q\in\cinfz{0,1}$ with
compact support. Again using Proposition \plref{S2.2} we deform
$q$ to $0$.

Thus it remains to prove the assertion for the operator $\-\frac{d^2}{dx^2}+z$
and $\nu_0, \nu_1\in \{\pm 1/2\}$. If $\nu_0=-1/2$ or
$\nu_1=-1/2$, in view of Proposition \plref{S2.3} it is enough to consider
the Neumann condition $f'(0)=0$ (resp. $f'(1)=0$). Repeating the
argument of the beginning of the proof we are left with the following
three operators:
\begin{itemize}
\item[1.] $D_1=-\frac{d}{dx^2}$ on $\{ f\in H^2[0,1]\,|\, f(0)=f(1)=0\}$,
\item[2.] $D_2=-\frac{d}{dx^2}$ on $\{ f\in H^2[0,1]\,|\, f(0)=0, f'(1)=0\}$,
\item[3.] $D_3=-\frac{d}{dx^2}+z$ on $\{ f\in H^2[0,1]\,|\, f'(0)=f'(1)=0\},
   z>0$.
\end{itemize}

1. $D_1=L_{1/2}$ and by Theorem \plref{S3.3} we have $\detz(D_1)=2$.
Moreover, $\varphi(x)=x, \psi(x)=1-x$, hence
   $$\frac{\pi}{2\Gamma(3/2)^2}W(\psi,\varphi)=2.$$

2. We have $\spec(D_2)=\{ (n+1/2)^2\pi^2\,|\, n\ge 0\}$, thus
$$\zeta_{D_2}(s)=\pi^{-2s}\sum_{n=1}^\infty (n+1/2)^{-2s}=
  \pi^{-2s}(2^{2s}-1)\zeta_R(2s),$$
and
$$\zeta'_{D_2}(0)=2\,\log 2\, \zeta_R(0)=-\log 2,$$
hence $\detz(D_2)=2.$ Furthermore, $\varphi(x)=x, \psi(x)=1$, thus
$$\frac{\pi\,W(\psi,\varphi)}{2^0 \Gamma(1/2)\Gamma(3/2)} =2.$$

3. Since the result is already proved for $D_1$ one finds
$$\detz(D_1+z)=2\frac{\sinh(\sqrt{z})}{\sqrt{z}}.$$
Furthermore, since $\spec(D_3)=\spec(D_1)\cup \{0\}$ we have
$$\detz(D_3+z)=2\sqrt{z} \sinh(\sqrt{z}).$$
On the other hand, we have $\varphi(x)=\cosh(\sqrt{z} x),
\psi(x)=\cosh(\sqrt{z}(x-1))$ and
$$W(\psi,\varphi)=\sqrt{z}\sinh(\sqrt{z})$$
and we are done.\hfill $\Box$\par\addvspace{0.25cm}

We single out the special case in which the Sturm--Liouville operator
can be factorized: let
\begin{equation}
   d:= \frac{d}{dx}+ S(x),
\end{equation}
where $S\in \cinf{(0,1)}$ such that
\begin{equation}
   S(x)= \frac{s_0}{x} +S_1(x)= \frac{s_1}{1-x}+ S_2(x)
\end{equation}
with $S_1\in\cinf{[0,1)}, S_2\in\cinf{(0,1]}$. Put
\begin{equation}
  l:= d^td=-\frac{d}{dx^2}+S^2-S'.
\end{equation}
Note that
\begin{equation}
  \nu_0=|s_0+1/2|,\quad \nu_1=|s_1-1/2|.
  \label{G4.3}
\end{equation}

Then the Friedrichs extension of $l$, $L:=l^\cf$, equals
$d^t_\max d_\min$ (cf. \cite[Lemma 3.1]{BL93}).

\begin{proposition}\label{S4.2} The $\zeta$--regularized determinant of
$L$ is given by the following formulas:

\noindent
$s_0\le -1/2, s_1\ge 1/2${\rm :}\quad $\detz(L)=0$,

\noindent
$s_0>-1/2, s_1<1/2${\rm :}
$$\detz(L)=\frac{\pi}{2^{\nu_0+\nu_1-2}\Gamma(\nu_0)\Gamma(\nu_1)}
  \exp\Big(-\regint_0^1 S(t) dt\Big)\;\int_0^1 \exp
     \Big(2\regint_0^x S(t) dt\Big) dx,$$

\noindent
$s_0>-1/2, s_1\ge 1/2${\rm :}
$$\detz(L)=\frac{\pi}{2^{\nu_0+\nu_1-1}\Gamma(\nu_0)\Gamma(\nu_1+1)}
   \exp\Big(\regint_0^1 S(t) dt\Big),$$

\noindent
$s_0\le -1/2, s_1<1/2${\rm :}
$$\detz(L)=\frac{\pi}{2^{\nu_0+\nu_1-1}\Gamma(\nu_0+1)\Gamma(\nu_1)}
   \exp\Big(-\regint_0^1 S(t) dt\Big).$$
\end{proposition}

\begin{proof} We put
\begin{equation}\begin{array}{rcl}
   \DST h(x)&:=&\DST\exp  \Big(-\regint_0^xS(t)dt\Big)= x^{-s_0}
            \exp\Big(-\int_0^x S_1(t) dt\Big).\\
       &=&\DST(1-x)^{s_1}  \exp\Big(-\regint_0^1 S(t) dt\Big)
           \exp\Big(\regint_x^1 S_2(t) dt\Big).
        \end{array}
\end{equation}
We have $dh=0$. Since $L=d_\max^t d_\min=d_\min^* d_\min$,
we have $\ker L=\ker d_\min$ and $\ker d_\min$ is non--trivial iff
\begin{eqnarray*}
       h(x)&=&o(x^{1/2} |\log x|^{1/2}),\, x\to 0,\\
       h(x)&=& o((1-x)^{1/2} |\log (1-x)|^{1/2}),\, x\to 1,
\end{eqnarray*}
(cf. \cite[Lemma 3.2]{BS88}), thus $\ker L\not=0$ iff $s_0\le -1/2$ and
$s_1\ge 1/2$.

$s_0>-1/2, s_1<1/2$: We put
\begin{eqnarray*}
   \varphi(x)&=& (2s_0+1) h(x) \int_0^x h(y)^{-2} dy,\\
   \psi(x)&=&  (1-2s_1) \exp \Big(-\regint_0^1 S(t) dt\Big)\; h(x) \int_x^1 h(y)^{-2} dy.
\end{eqnarray*}
It is easy to check that $\varphi$ is normalized at $0$ and $\psi$ is
normalized at $1$ and
$$W(\psi,\varphi)=(2s_0+1)(2s_1+1) \exp\Big(-\regint_0^1 S(t) dt\Big)
   \int_0^1 h(y)^{-2} dy.$$
Using Theorem \plref{S4.1} we reach the conclusion.

$s_0>-1/2, s_1\ge 1/2$: We put
\begin{eqnarray*}
   \varphi(x)&=& (2 s_0+1) h(x) \int_0^x h(y)^{-2} dy,\\
   \psi(x)&=&  \exp \Big(\regint_0^1 S(t) dt\Big)\; h(x).
\end{eqnarray*}
Then $\varphi$ is normalized at $0$ and $\psi$ is normalized at $1$ and
$$W(\psi,\varphi)=(2s_0+1) \exp \Big(\regint_0^1 S(t) dt\Big) $$
and again we reach the conclusion using Theorem \plref{S4.1}.

$s_0\le -1/2, s_1<1/2$: \\This is proved analogously to the case
$s_0>-1/2, s_1\ge 1/2.$\end{proof}

As a classical example we discuss

\newsubsection{The Jacobi differential operator}

For $\ga, \gb>-1$, the Jacobi polynomials, $P_n^{(\ga,\gb)}, n\ge 0$,
form a complete orthogonal set in the Hilbert space
\begin{equation}
  \ch:= L^2([-1,1], (1-x)^\ga(1+x)^\gb dx).
\end{equation}
We put
\begin{equation}
\varrho(x):=(1-x)^\ga (1+x)^\gb,\quad p(x):=1-x^2.
\end{equation}

$P_n^{(\ga,\gb)}$ satisfies the differential equation \cite[p. 258]{R}
\begin{equation}
  -\frac{1}{\varrho(x)} \frac{d}{dx} \varrho(x)p(x) \frac{d}{dx}
   P_n^{(\ga,\gb)}(x)=n(n+\ga+\gb+1) P_n^{(\ga,\gb)},
\end{equation}
thus $P_n^{(\ga,\gb)}$ are eigenfunctions of the operator
\begin{equation}
  j:= -\frac{1}{\varrho(x)} \frac{d}{dx} \varrho(x)p(x) \frac{d}{dx}
\end{equation}
and it is not difficult to see that the $P_n^{(\ga,\gb)}$ are in the
domain of the self--adjoint extension $J=d_\max^*d_\max$, where
$d=\sqrt{\varrho}\frac{d}{dx}$. Hence we have
\begin{equation}
\spec(J)=\{\, n(n+\ga+\gb+1)\,|\, n=0,1,2,\ldots\}.
  \label{G4.4}
\end{equation}

Next we transform $j$ into an operator in $L^2[0,1]$. We put
\begin{equation}
\kappa(x):=\frac{1}{\pi} \arcsin(x)+1/2,\quad -1\le x\le 1.
\end{equation}

Now a straightforward calculation shows:

\begin{lemma}\label{S4.3} The map
$$\Phi:L^2[0,1]\longrightarrow L^2([-1,1],\varrho(x) dx),
   f\mapsto \sqrt{\kappa'} \varrho^{-1/2} f\circ \kappa$$
is unitary and we have
$$\Phi^* j \Phi =\frac{1}{\pi^2} d^t d$$
with
\begin{equation}
    d= \frac{d}{dx} -\frac{\pi}{2}(1+\ga+\gb)\cot(\pi x)+\frac{\pi}{2}
         \frac{\ga-\gb}{\sin(\pi x)}.
\end{equation}
\end{lemma}

Thus we are almost in the situation of Proposition \plref{S4.2}, except
that $0$ is an eigenvalue of $J$ and hence of $\Phi^* J \Phi$. Note
that
$$\Phi^* J \Phi=\frac{1}{\pi^2} d_\max^* d_\max,$$
hence
\begin{equation}
      \spec(\Phi^* J \Phi)=\spec(\frac{1}{\pi^2} d_\max d_\max^*)\cup
      \{0\}
      \label{G4.1}
\end{equation}
and $d_\max d_\max^*=(d^t_\min)^*d^t_\min=(d d^t)^\cf$. Now we have
\begin{eqnarray*}
    -d^t&=& \frac{d}{dx}+\frac{\pi}{2}(1+\ga+\gb)\cot(\pi x)-\frac{\pi}{2}
         \frac{\ga-\gb}{\sin(\pi x)}\\
        &=:& \frac{d}{dx} + S(x)=: d_{\ga,\gb}.
\end{eqnarray*}
and
\begin{equation}
  S(x)\sim \frac{1/2+\gb}{x},\, x\to 0,\quad S(x)\sim \frac{-1/2-\ga}{1-x},
\, x\to 1.
  \label{G4.5}
\end{equation}
We calculate $d_{\ga,\gb}^td_{\ga,\gb}$ explicitly:
\begin{eqnarray}
   d_{\ga,\gb}^td_{\ga,\gb}&=& -\frac{d^2}{dx^2}+\frac{\pi^2}{4}\Bigg(
     ((\ga+\gb+2)^2-1) \cot^2 \pi x + \frac{(\ga-\gb)^2}{\sin^2\pi x}
        \nonumber\\
    &&-2 (\ga-\gb)(\ga+\gb+2)\frac{\cos \pi x}{\sin^2 \pi x} +2(\ga+\gb+1)
       \Bigg)\label{G4.6}\\
     &=:& l_{\ga,\gb}.
\end{eqnarray}
Let $L_{\ga,\gb}:=l_{\ga,\gb}^\cf$ be the Friedrichs extension of $l_{\ga,\gb}$.
$L_{\ga,\gb}$ obviously makes sense even for $\ga=-1$ or $\gb=-1$.

\begin{lemma}\label{S4.7} We have
\begin{equation}
    \spec(L_{\ga,\gb})=\{ \pi^2\, n(n+\ga+\gb+1)\,|\, n=1,2,\ldots\},
    \quad \ga,\gb\ge -1.
    \label{G4.2}
\end{equation}
\end{lemma}
\begin{proof} By Proposition \plref{S4.2} $L_{\ga,\gb}$ is invertible
for $(\ga,\gb)\not=(-1,-1)$. Hence, for $\ga,\gb>-1$, the assertion
follows from \plref{G4.4} and \plref{G4.1}.

Now, a straightforward calculation shows:
\begin{eqnarray}
       d_{\ga,\gb}d_{\ga,\gb}^t&=& -\frac{d^2}{dx^2}+\frac{\pi^2}{4}\Bigg(
     ((\ga+\gb)^2-1) \cot^2 \pi x + \frac{(\ga-\gb)^2}{\sin^2\pi x}
        \nonumber\\
    &&-2 (\ga-\gb)(\ga+\gb)\frac{\cos \pi x}{\sin^2 \pi x} -2(\ga+\gb+1)
       \Bigg)\nonumber\\
    &=& l_{\ga-1,\gb-1}-\pi^2(\ga+\gb).\label{G4.8}
\end{eqnarray}
For $\ga\ge 0, \gb\ge 0$ we infer from \cite[Lemma 3.2]{BS88} and
\myref{G4.5} that $d_{\ga,\gb,\max}=d_{\ga,\gb,\min}$ and hence we find
\begin{eqnarray}
   L_{\ga-1,\gb-1}-\pi^2(\ga+\gb)&=&(d_{\ga,\gb}d_{\ga,\gb}^t)^\cf
      =d_{\ga,\gb,\max}(d_{\ga,\gb}^t)_\min\label{G4.9}\\
      &=&d_{\ga,\gb,\min}(d_{\ga,\gb,\min})^*.\nonumber
\end{eqnarray}
Moreover, from Proposition \plref{S4.2} we infer
$\ker (d_{\ga,\gb}^t)_\min\not=0$, hence
$0\in \spec(L_{\ga-1,\gb-1}-\pi^2(\ga+\gb))$ and \myref{G4.9}
implies
\begin{eqnarray*}
   \spec(L_{\ga-1,\gb-1}-\pi^2(\ga+\gb))&=&
     \spec((d_{\ga,\gb,\min})^*d_{\ga,\gb,\min})\cup \{0\}\\
     &=&\spec(L_{\ga,\gb})\cup \{0\},
\end{eqnarray*}
thus
\begin{eqnarray*}
   \spec(L_{\ga-1,\gb-1})&=&\{\pi^2 n(n+\ga+\gb+1)+\pi^2(\ga+\gb)\,|\,
    n=0,1,2,\ldots\}\\
    &=&\{\pi^2 n(n+ (\ga-1)+(\gb-1)+1)\,|\, n=1,2,\ldots\}.
\end{eqnarray*}
\end{proof}

We calculate $\detz(L)$ using Proposition \plref{S4.2}:
\begin{small}
$$\renewcommand{\arraystretch}{2.5}\begin{array}{l}
  \DST-\log h(x)=\regint_0^x S(t) dt=\LIM_{\eps\to 0} \int_\eps^x S(t) dt\\
   \DST\quad=\LIM_{\eps\to 0} \left\{ \frac{1+\ga+\gb}{2}\log \sin(\pi t)
      \Big|_{t=\eps}^{t=x}+\frac{\gb-\ga}{2} \log \Big(\frac{1-\cos(\pi t)}%
         {\sin(\pi t)}\Big)\Big|_{t=\eps}^{t=x}\right\} \\
   \DST\quad= \frac{1+\ga+\gb}{2}\log \sin(\pi x)- \frac{1+\ga+\gb}{2} \log\pi+
   \frac{\gb-\ga}{2} \log \Big(\frac{1-\cos(\pi x)}{\sin(\pi x)}\Big)
   +\frac{\ga-\gb}{2}\log\frac{\pi}{2}, \\
   \DST\quad =(1/2+\gb)\log\sin(\frac \pi 2 x)+(1/2+\ga)\log\cos(\frac \pi 2 x)-
      (1/2+\gb)\log \pi +(1/2+\gb)\log 2,
    \end{array}$$
\end{small}
and
$$\regint_0^1 S(t) dt=
   \LIM_{x\to 1} \regint_0^x S(t) dt= (\ga-\gb) \log \frac{\pi}{2}.$$
Moreover
\begin{eqnarray*}
   \int_0^1 h(x)^{-2} dx&=& \pi^{-1-2\gb} 2^{1+2\gb} \int_0^1
      \big(\sin(\frac \pi 2 x)\big)^{1+2\gb} \big(\cos \frac \pi 2 x\big)^{1+2\ga} dx\\
      &=&\pi^{-2-2\gb} 2^{1+2\gb} \frac{\Gamma(\ga+1)\Gamma(\gb+1)}{\Gamma(2+\ga+\gb)}.
\end{eqnarray*}

Using Proposition \plref{S4.2} we  have proved:

\begin{proposition}\label{S4.4} For $\ga,\gb\ge -1,$ the determinant of the Jacobi--operator
$L_{\ga,\gb}$ is
$$\detz(L_{\ga,\gb})=\frac{2 \pi^{-1-\ga-\gb}}{\Gamma(2+\ga+\gb)}.$$
\end{proposition}

Note that if $\ga=\gb=-1$ then $\detz(L_{\ga,\gb})=0$. Since $\Gamma$ has
a pole at $0$, the formula also covers this case.

Since we know $\spec(L)$ explicitly, this result can also be proved directly.
This \linebreak in fact leads to an alternative proof of Theorem \plref{S4.1}, that
does not make use of section \sec1.

For doing this, we introduce the function
\begin{equation}
\zeta_\lambda(s):= \sum_{n=1}^\infty n^{-s}(n+\gl)^{-s},\quad
  \gl>-1,\; \Re s>1/2.
\end{equation}

\begin{lemma}\label{S4.5} $\zeta_\gl$ has a meromorphic continuation to $\C$.
$\zeta_\gl$ is regular at $s=0$ and we have
$$\zeta_\gl(0)=-\frac{1+\gl}{2},\quad \zeta_\gl'(0)=-\log {2\pi}+
  \log \Gamma(\gl+1).$$
\end{lemma}
\begin{proof} That $\zeta_\gl$ has a meromorphic continuation is well--known.
A simple way of seeing this is
\begin{eqnarray*}
    \zeta_\gl(s)&=& \sum_{1\le n\le \gl+1} n^{-s} (n+\gl)^{-s}+
       \sum_{n>\gl+1} n^{-2s}(1+\frac{\gl}{n})^{-s}\\
       &=& \sum_{1\le n\le \gl+1} n^{-s} (n+\gl)^{-s}+
         \sum_{k=0}^\infty {-s \choose k} \sum_{n>\gl+1} n^{-2s-k} \gl^k
\end{eqnarray*}
and the right hand side is a meromorphic function in the whole plane.
Moreover we have
$$\zeta_\gl(s)=\zeta_R(2s)-\gl s \zeta_R(2s+1) +
    \sum_{n=1}^\infty \Big[n^{-s}(n+\gl)^{-s}-n^{-2s}+\gl s n^{-2s-1}\Big],$$
Since
$$|n^{-s}(n+\gl)^{-s} - n^{-2s} +\gl s n^{-2s-1}|=O(s \gl n^{-2 \Re s -2}),$$
this shows that $\zeta_\gl$ is regular at $s=0$. We find
$$\zeta_\gl(0)=\zeta_R(0)-\frac \gl 2 \Res_{s=1} \zeta_R(s)=-\frac{1+\gl}{2}$$
and
\begin{eqnarray*}
   \zeta'_\gl(0)&=&2\zeta_R'(0)-\gl \frac{d}{ds}\Big|_{s=0}
       \Big( s\zeta_R(2s+1)\Big)
     -\sum_{n=1}^\infty \Big[\log(n+\gl)-\log n-\frac{\gl}{n}\Big]\\
     &=&- \log (2 \pi)-\gl \gamma-\sum_{n=1}^\infty
       \Big[ \log(1+\frac{\gl}{n})-
       \frac \gl n \Big]\\
     &=&-\log(2\pi)+\log \Gamma(\gl)+\log \gl\\
     &=&-\log(2\pi)+\log \Gamma(\gl+1).
\end{eqnarray*}
\end{proof}

In view of \plref{G4.2} we have
$$\zeta_L(s)=\pi^{-2s} \zeta_{1+\ga+\gb}(s),$$
thus
\begin{eqnarray*}
   \log\detz(L)&=& 2\,\log \pi\, \zeta_{1+\ga+\gb}(0)-\zeta_{1+\ga+\gb}'(0)=
     \log\frac{2 \pi^{-1-\ga-\gb}}{\Gamma(2+\ga+\gb)}.
\end{eqnarray*}

As promised, we sketch a

\medskip

\noindent
{\bf Second proof of Theorem \ref{S4.1}.}\quad
\enspace
As in the first proof, use Propositions \plref{S2.2} and \plref{S2.3}
to show that
$$\detz(L)=c(\nu_0,\nu_1) W(\psi,\varphi)$$
with some constant $c(\nu_0,\nu_1)$ depending only on $\nu_0, \nu_1$.
Since Proposition \plref{S4.4} can be proved directly, we may use it
to show that
$$c(\nu_0,\nu_1)=\frac{\pi }{2^{\nu_0+\nu_1}\Gamma(\nu_0+1)\Gamma(\nu_1+1)},
\quad\mbox{\rm if}\quad \nu_0,\nu_1\ge 0,\; \nu_0+\nu_1>0.$$

For fixed $\nu_0$ choose a symmetric potential $q(x)=q(1-x)$ with
$\nu_0=\nu_0(q)$. Then an easy calculation shows that
the eigenvalues of $(-\frac{d^2}{dx^2}+q)^\cf$ consist of the union
of the eigenvalues of  $-\frac{d^2}{dx^2}+q$ on $[0,1/2]$ with Dirichlet
and Neumann boundary conditions at $1/2$ and the original boundary condition
at $0$. A direct calculation shows that this implies
$$c(\nu_0,\nu_0)=\frac 12 c(\nu_0,1/2) c(\nu_0,-1/2),$$
proving
\begin{eqnarray*}
    c(\nu_0,-1/2)&=& 2 \frac{c(\nu_0,\nu_0)}{c(\nu_0,1/2)}\\
       &=& \frac{\pi}{2^{\nu_0-1/2} \Gamma(\nu_0+1)\Gamma(1/2)}
\end{eqnarray*}
and
$$c(1/2,-1/2)=2.$$
$c(-1/2,-1/2)$ is now calculated as in the first proof.
The case $\nu_0=\nu_1=0$ has to be treated separately. We leave the
details to the reader.\hfill $\Box$\par\addvspace{0.25cm}

\newsubsection{$\detz(L+z)$ as an infinite product}

\begin{proposition}\label{S4.6} Let $L$ be a semibounded invertible self--adjoint
operator in the Hilbert space $\ch$ satisfying \myref{G1.1} and
\myref{G1.2}. Let $(\gl_n)_{n\in\N}$ be the eigenvalues of $L$.
Then we have for $z\in\C$
$$\detz(L+z)=\detz(L)\prod_{n=1}^\infty (1+\frac{z}{\gl_n}).$$
\end{proposition}

\begin{proof} In view of \myref{G1.1} left and right hand side of the
equation are entire holomorphic functions and we find
\begin{eqnarray*}
    \frac{d}{dz} \log \detz(L+z)&=& \Tr((L+z)^{-1})
           =\sum_{n=1}^\infty (\gl_n+z)^{-1}\\
       &=&\sum_{n=1}^\infty \frac{d}{dz} \log(1+\frac{z}{\gl_n})
       =\frac{d}{dz} \log\,\prod_{n=1}^\infty (1+\frac z{\gl_n}).
\end{eqnarray*}
Since the assertion is obviously true for $z=0$
we reach the conclusion.\end{proof}

We apply this formula to a regular singular Sturm--Liouville operator
$L=(l,R_0,R_1)$. Let $\varphi(\cdot,z), \psi(\cdot,z)$ be the
normalized solutions for $L+z^2$. Then applying Theorem \plref{S4.1} and
the preceding proposition we find
$$   \prod_{n=1}^\infty (1+\frac{z^2}{\gl_n})=\frac{\detz(L+z^2)}{\detz(L)}=
   \frac{W(\psi(\cdot,z),\varphi(\cdot,z))}{W(\psi(\cdot,0),\varphi(\cdot,0))}.$$
In the case of the operator $-\frac{d^2}{dz^2}$ with Dirichlet boundary 
conditions,
this is essentially the product expansion of $\sinh$, namely we have
$\varphi(x,z)=\sinh(xz)/z$, thus
\begin{equation}
\prod_{n=1}^\infty (1+\frac{z^2}{n^2 \pi^2})=\frac{\varphi(1,z)}{\varphi(1,0)}=
    \frac{\sinh(z)}{z}.
\end{equation}

More generally, let $(\gl_{n,\nu})_{n\in\N}$ be the zeros of the
Bessel function $J_\nu$. Then we have
$$\spec(L_\nu)=\{  \gl_{n,\nu}^2\,|\, n\in \N\}.$$
Furthermore
$$\varphi(x,z)=2^\nu \Gamma(\nu+1) z^{-\nu} \Inu(xz), \quad
   \varphi(x,0)=x^{\nu+1/2},$$
thus
$$\prod_{n=1}^\infty (1+\frac{z^2}{\gl_{n,\nu}^2})=
    \frac{\varphi(1,z)}{\varphi(1,0)}=2^\nu z^{-\nu} \Gamma(\nu+1)\Inu(z)$$
or
\begin{equation}
\Inu(z)=\frac{(z/2)^\nu}{\Gamma(\nu+1)}\prod_{n=1}^\infty(1+
   \frac{z^2}{\gl_{n,\nu}^2}).
\end{equation}
Of course, this formula is classical \cite[Sec. 15.41 (3)]{Watson}.

\newsection{An open problem}

\def\sec5{5}

We briefly outline our initial motivation for proving
Theorem \plref{S4.1}.

Let $M^m$ be a compact Riemannian manifold. Then the analytic torsion
\cite{RS} is defined by
\begin{equation}
\log T(M)=\frac 12 \sum_{i=0}^m (-1)^i i \zeta_i'(0),
\end{equation}
where $\zeta_i(s)$ denotes the $\zeta$--function of the Laplacian on
$i$--forms.

The celebrated Cheeger--M\"uller theorem \cite{C,Muller} identifies $T(M)$ with
a purely combinatorial object, the combinatorial torsion of $M$.

\begin{problem}
Is there an analogue of the Cheeger--M\"uller
theorem for a suitable class of pseudomanifolds?
\end{problem}

Few attempts have been made in this direction. A. Dar \cite{Dar1,Dar2}
defined and investigated Reidemeister torsion for intersection cohomology
and one might expect that the intersection cohomology should show up on
the combinatorial side. On the analytic side, only manifolds with cone--like
singularities have been considered. Using work of Cheeger, A. Dar
proved

\begin{proposition}\label{S5.1} {\rm \cite{Dar1}}\quad Let $M$ be a compact manifold
with cone--like singularities. Then $T(M)$ exists.
\end{proposition}

The only thing one has to show is that the meromorphic function
\begin{equation}
\sum_{i=0}^m (-1)^i i \zeta_i(s)
\end{equation}
has no pole at $0$. A priori this function has a simple pole at $0$
due to a $\log$--term in the heat asymptotics. However, the sum
\begin{equation}
\sum_{i=0}^m (-1)^i i \Res_{s=0} \zeta_i(s)
\end{equation}
turns out to be $0$.

An interesting approach to the Cheeger--M\"uller theorem is the recent
work of Vishik \cite{Vishik}, who proves a gluing formula for
analytic torsion norms. Adopting this approach, for proving an
analogue of the Cheeger--M\"uller theorem for manifolds with
cone--like singularities it would be enough to compare
analytic and the (hypothetical) R--torsion for the model cone
$C(N)$ over a compact manifold $N$. At least this would indicate what
a result could look like.

More precisely, let
$$C(N)=(0,1)\times N$$
be the model cone over $N$ with metric
$$g= dx^2\moplus x^2 g_N.$$

On the face $\{1\}\times N$ we impose relative boundary conditions. Then
separation of variables shows that the Laplacian on $i$--forms is
an infinite sum of operators
$$L_\gl=-\frac{d^2}{dx^2}+\frac{\gl}{x^2},$$
where the $\gl$'s are essentially the eigenvalues of the Laplacian on $N$.
This is the reason why the author considered Theorem \plref{S4.1}.
However, for calculating $T(C(N))$ it is not enough to know
$\detz(L_\gl)$, since one has to deal with an infinite sum of
operators. We leave this as a problem

\begin{problem} Calculate $T(C(N))$ for relative/absolute
boundary conditions. 
\end{problem}

Together with Vishik's result, the solution to this problem should
lead to a Cheeger--M\"uller type result for manifolds with cone--like
singularities.

\begin{acknowledgements}
This work was written while I was visiting the Ohio State University,
Columbus, Ohio.

\begin{sloppypar}
It is my pleasure to thank the people at the Mathematics Department
of OSU for making my stay so fruitful and enjoyable. In particular I
would like to thank {\sc Henri Moscovici} and {\sc Bob Stanton} for their help
and encouragement.

This work is a byproduct of a somewhat frustrating project which
{\sc Henri Moscovici} and I had started and with which
we became stuck at the very beginning. Nevertheless this paper would not
have been written without Henri's vision of a Cheeger--M\"uller
Theorem for pseudomanifolds (cf. Sec. {\rm \secfive}). I wish to
thank {\sc Dan Burghelea} for showing me the preprint {\rm \cite{BFK2}} and
{\sc Thomas Kappeler} for many useful discussions.
Moreover, I owe the references {\rm \cite{Dar1,Dar2}} to {\sc Jeff Cheeger}
for which I would like to thank him.

Finally, I express my gratitude to the anonymous referee
for helpful comments and for pointing out an error to me.
\end{sloppypar}
\end{acknowledgements}

\address{Humboldt--Universit\"at zu Berlin\\
Institut f\"ur Mathematik\\
Unter den Linden 6\\
D - 10099 Berlin\\
Germany\\
e-mail: lesch@mathematik.hu-berlin.de}
\address{ }

\end{document}